\documentclass[journal]{IEEEtran}

\usepackage[T1]{fontenc}
\usepackage[utf8]{inputenc}
\usepackage{amsmath}
\usepackage{textcomp}
\usepackage{upgreek}
\usepackage{hyperref}
\usepackage{cleveref}
\usepackage{commath}
\usepackage{comment}
\usepackage{float}
\usepackage{graphicx}
\usepackage{listings}
\usepackage{tikz}
\usepackage[makeroom]{cancel}
\usepackage{gensymb}
\usepackage{makecell}
\usepackage{booktabs}
\usepackage{dcolumn}
\usepackage{wrapfig}
\usepackage{mathtools}
\usepackage{siunitx}
\usepackage[caption=false,font=footnotesize,position=top]{subfig}
\usepackage[official]{eurosym}

\newcommand{\N}[1]{\mathcal{N}_{#1}}
\newcommand{\M}{\mathcal{M}}
\newcommand{\LO}{\mathcal{L}}
\newcommand{\I}[2]{\mathcal{I}_{#1}^{#2}}

\newcommand{\T}{\mathcal{T}}
\newcommand{\SC}{\mathcal{S}}
\newcommand{\X}{\mathcal{X}}
\newcommand{\Y}{\mathcal{Y}}
\renewcommand{\vec}[1]{\boldsymbol{\mathbf{#1}}}
\usepackage{cite}

\ifCLASSINFOpdf
\else
\fi

\usepackage{nomencl}
\makenomenclature

\usepackage{etoolbox}
\renewcommand\nomgroup[1]{%
  \item[\bfseries
  \ifstrequal{#1}{S}{Sets}{%
  \ifstrequal{#1}{V}{Variables}{%
  \ifstrequal{#1}{P}{Parameters}{}}}%
]}

\begin{document}
\bstctlcite{MyBSTcontrol}

%
\title{Procurement of Spinning Reserve Capacity in a Hydropower Dominated System Through Mixed Stochastic-Robust Optimization}
%
%
%

\author{Christian~Øyn~Naversen,~\IEEEmembership{Student Member,~IEEE,}
        Hossein~Farahmand,~\IEEEmembership{Member,~IEEE,}
		Arild~Helseth,~\IEEEmembership{Member,~IEEE}
\thanks{ C. Ø. Naversen and H. Farahmand are with the Department of Electric Power Engineering, Norwegian University of Science and Technology, 7491 Trondheim, Norway. E-mail: christian.naversen@ntnu.no, hossein.farahmand@ntnu.no}%
\thanks{A. Helseth is with the Department of Energy Systems, SINTEF Energy Research, Sem Sælands vei 11, 7034 Trondheim, Norway. E-mail: Arild.Helseth@sintef.no}}
\maketitle

\begin{abstract}
As the penetration of variable renewable power generation increases in power systems around the world, system security is challenged. Ensuring that enough reserve capacity is available to balance the increased forecast errors introduced by wind and solar power are and will be an important challenge to solve for system operators. This is also true in systems that already have flexible balancing resources, such as hydropower. The main challenge in this case will be the coordination of the reserve procurement between connected hydropower plants, as the water flow is changed when the reserve capacity is activated. However, the activation step is often ignored or simplified in hydropower scheduling models that include reserve capacity procurement. In this paper, a two-stage model for scheduling power and procuring symmetric spinning reserves in a hydropower system with uncertain net load is proposed to capture the effect of geographic reserve capacity coordination. To model the uncertainty in the net load, a new mixed stochastic-robust optimization model is proposed to achieve cost efficiency and system security. We show that this proposed model outperforms its natural contenders in the given case study, and so does not suffer from the typical overly conservative nature of pure robust models.
\end{abstract}



%
\IEEEpeerreviewmaketitle

%
%
%
%

\nomenclature[S,01]{$\mathcal{M}$}{Hydropower modules}
\nomenclature[S,02]{$\mathcal{N}_{m}$}{Discharge segments in module $m$}
\nomenclature[S,03]{$\mathcal{I}^{d/b/o}_{m}$}{Modules that discharge/bypass/spill water into module $m$}
\nomenclature[S,04]{$\mathcal{O}^{d/b/o}_{m}$}{Modules that module $m$ discharges/bypasses/spills water into}
\nomenclature[S,05]{$\mathcal{T}$}{Time periods}
\nomenclature[S,06]{$\mathcal{S}$}{Balancing scenarios}
\nomenclature[S,07]{$\mathcal{J}$}{Robust balancing scenarios}
\nomenclature[S,08]{$\X$}{First-stage feasibility constraints}
\nomenclature[S,09]{$\Y$}{Balancing-stage feasibility constraints}
\nomenclature[S,10]{$\Omega$}{Dual feasibility constraints for the balancing stage}
\nomenclature[S,11]{$\LO$}{Robust uncertainty set}

\nomenclature[V,01]{$q^{d}_{mnt}$}{Discharge through segment $n$ in module $m$ during period $t$}
\nomenclature[V,02]{$q^{b}_{mt}$}{Bypass from module $m$ during period $t$}
\nomenclature[V,03]{$q^{o}_{mt}$}{Overflow from module $m$ during period $t$}
\nomenclature[V,04]{$q^{in}_{mt}$}{Water flow into module $m$ during period $t$}
\nomenclature[V,05]{$q^{out}_{mt}$}{Water flow out of module $m$ during period $t$}
\nomenclature[V,06]{$v_{mt}$}{Volume in module $m$ at the beginning of period $t$}
\nomenclature[V,07]{$v_{m,T+1}$}{End volume in module $m$}
\nomenclature[V,08]{$p_{mt}$}{Power generated in module $m$ during period $t$}
\nomenclature[V,09]{$r_{mt}$}{Symmetric spinning capacity reserved on module $m$ during period $t$}
\nomenclature[V,10]{$s_{t}^{\pm}$}{Load shed and power spill during period $t$}
\nomenclature[V,11]{$u_{t}^{\pm}$}{Load deviation in upward and downward direction during period $t$}
\nomenclature[V,12]{$\vec{x}$}{Vector of all first-stage variables}
\nomenclature[V,13]{$\vec{y}$}{Vector of all balancing-stage variables}
\nomenclature[V,14]{$\vec{\phi}$}{Vector of all dual balancing-stage variables}

\nomenclature[P,01]{$WV_{m}$}{End value of water in module $m$}
\nomenclature[P,02]{$C^{b/o}$}{Penalty cost for bypassing/spilling water}
\nomenclature[P,03]{$C^{\pm}$}{Penalty cost for shedding load/spilling power}
\nomenclature[P,04]{$Q^{d}_{mn}$}{Maximal discharge through segment $n$ in module $m$}
\nomenclature[P,05]{$Q^{b/o}_{m}$}{Maximal flow through bypass/spill gate from module $m$}
\nomenclature[P,06]{$I_{mt}$}{Inflow to module $m$ during period $t$}
\nomenclature[P,07]{$V^{0}_{m}$}{Initial storage level of module $m$}
\nomenclature[P,08]{$V_{m}$}{Maximal storage capacity of module $m$}
\nomenclature[P,09]{$E_{mn}$}{Energy conversion factor for segment $n$ in module $m$}
\nomenclature[P,10]{$P_{m}$}{Maximal production capacity of module $m$}
\nomenclature[P,11]{$L_{t}$}{Forecasted system net load during period $t$}
\nomenclature[P,12]{$R_{t}$}{Spinning system reserve requirement during period $t$}
\nomenclature[P,13]{$F_{t}$}{Length of period $t$}
\nomenclature[P,14]{$T$}{Number of periods in $\T$}
\nomenclature[P,15]{$\Gamma$}{Budget of uncertainty}
\nomenclature[P,16]{$\Lambda$}{Maximal net load deviation}

\section{Introduction}

\IEEEPARstart{H}{ydropower} is a valuable balancing asset for any power system, as it is flexible and cheap compared to thermal generation technologies. As the share of variable renewable energy sources in power systems across the world increases, so does the need for balancing capacity and energy. Although hydropower is well suited to help regulate the system, the technical constraints and cascaded topology must be considered to realistically estimate this regulating capability. The watercourse connects all hydropower units in space and time, and so the balancing actions of a single unit will impact the whole system. This is a challenge when considering spinning reserve capacity allocation, as sufficient water must be available for production when the reserve capacity is activated. Additionally, enough storage capacity in the reservoirs below is essential to avoid spillage and the loss of potential energy. Another complicating aspect is the implicitly defined marginal cost of operating a hydropower plant. The stored water in each reservoir has an associated opportunity cost or water value, which in general depends on the state of the entire system. This makes the cost of procuring reserve capacity on a specific hydropower plant dependent on the capacity procured on the surrounding units.

Reserve capacity procurement and system balancing have been incorporated into hydropower scheduling models for cascaded systems in several ways. These features can be found in both long-term planning models \cite{Jaehnert2012, Abgottspon2014, Helseth2016, Hjelmeland2016,Lohndorf2013} and short-term operational models \cite{Kong2017,Kazempour2008, Chazarra2016,Gu2014}. The fundamental model in \cite{Jaehnert2012} sequentially clears the day-ahead market, reserve procurement and system balancing steps for Northern Europe. The activation of reserve capacity is based on the marginal cost of the hydropower units in their day-ahead position, but does not include hydrological constraints. The methods presented in \cite{Abgottspon2014} and \cite{Helseth2016} consider a producer participating in day-ahead energy and spinning reserve capacity markets under uncertainty in inflow and market prices within modified stochastic dual dynamic programming (SDDP) frameworks. Both works ensure that enough water is stored in the reservoir to produce the allocated reserve capacity, though activation is not directly modelled. In \cite{Hjelmeland2016} it is investigated how wind power can contribute to the provision of rotating reserves in a hydropower-dominated system by using the SDDP algorithm, but without considering reserve activation. The deterministic model presented in \cite{Kong2017} has a high degree of physical detail, and can model the reservation of all the different reserve capacity products in Norway. The total amount of reserve capacity to allocate in the system is exogenously given to the model, and is distributed among the hydropower units while optimizing the day-ahead market position. The probability of activation in the balancing markets modifies the expected income in the deterministic model in \cite{Kazempour2008}, and the work in \cite{Chazarra2016} is based on the assumption that a certain percentage of the reserve capacity sold to the market is activated by the system operator in every inflow and price scenario. The works in \cite{Lohndorf2013} and \cite{Gu2014} do not explicitly model the reserve capacity procurement, but considers system balancing through bidding into the day-ahead, intraday and real-time energy markets. To the best of the authors' knowledge, the importance of representing the activation of procured reserves in the energy and reserve scheduling for a large-scale hydropower system has not been addressed in detail in the literature, and so emerges as a gap in the existing research.

The medium-term hydrothermal model presented in \cite{Street2017} procures reserve capacity to ensure system security in the face of a $N-k$ security criterion. This is done by incorporating robust optimization into the SDDP framework, and marks one of several interesting optimization problems that employs both stochastic and robust optimization in an effort to reap the benefits of both methods. In contrast to stochastic optimization based on the expected objective value over a set of scenarios, robust optimization hedges the solution against the worst case realization of the uncertainty. Robust optimization has been widely and successfully applied to power system planning and operation problems in recent years \cite{Morales2014}. A large portion of the published scientific material has been related to the unit commitment problem under uncertainty, where the goal typically is to commit a sufficient number of thermal units to be able to balance real-time deviations \cite{VanAckooij2018}. These types of models are usually formulated as two-stage models \cite{Jiang2012,Bertsimas2013,Zhao2013,An2015,Blanco2017}, though single-stage \cite{Street2011} and multistage models \cite{Lorca2016,Lorca2017} also exist. Robust optimization has also been used in the context of hydropower scheduling under uncertainty, as in \cite{Street2017} and \cite{Apostolopoulou2018}. The combination of robust and stochastic optimization has been proposed in different ways. Stochastic and robust optimization may handle separate sources of uncertainty, such as generator availability and power prices in \cite{Dehghan2016} and variable power generation and power prices in \cite{Liu2016}. It is also possible to make hybrid models by taking a stochastic or a robust model and introducing some characteristics from the other approach. The work in \cite{Blanco2017} partitions the scenarios in a stochastic model into bundles where robust optimization is enforced within each bundle, while \cite{An2015} introduces several robust uncertainty sets to a robust model by weighting them in the objective function akin to scenario probabilities.

The model presented in this paper draws its inspiration from the unified stochastic-robust model presented in \cite{Zhao2013}. The same source of uncertainty is modelled by both stochastic and robust optimization in \cite{Zhao2013} by introducing a weight $\beta$ of the average scenario cost and $1-\beta$ of the robust worst case cost in the objective function. This represents a direct integration of both the stochastic and the robust optimization methods in a single problem, which opens up new and interesting avenues of hybridization. We define a new variant of this model and apply it to the hydropower scheduling and reserve procurement problem under net load uncertainty.  We show that a good choice of $\beta$ leads to a model that is both more robust and less costly compared to its deterministic, pure robust, pure stochastic and unified stochastic-robust counterparts. The new model is less complex than the original unified model, as it separates out the calculation of the robust elements. This is done by leveraging the popular column-and-constraint generation (CCG) solution technique (see \cite{Zeng2013,LongZhao2012}) as a scenario generator, which provides robust scenarios for the proposed mixed stochastic-robust model. A detailed look at the impact of the uncertainty modelling on the hydropower reserve capacity procurement is also provided. In short, the contributions of this paper are considered twofold: 

\begin{enumerate}
    \item A new mixed stochastic-robust optimization model is proposed. This is shown to be an improvement over the existing model formulations when applied to the hydropower scheduling and reserve procurement problem under uncertainty.
    \item The importance of coordinating the reserve capacity procurement in a cascaded hydropower system is demonstrated with a realistic example based on a Norwegian watercourse, where the differences in performance between deterministic and mixed stochastic-robust optimization are analyzed.
\end{enumerate}

The rest of the paper is organized into three parts: \Cref{model_section} details the modelling of the optimization problem formulations, a case study is presented in \Cref{case_study_section}, and concluding remarks are found in \Cref{conclusion_section}. \Cref{model_section} is split into subsections describing the deterministic day-ahead scheduling problem (\Cref{deterministic_section}), the system balancing problem (\Cref{bal_section}), the pure stochastic and pure robust two-stage problems (\Cref{two_stage_section}), and the unified and new mixed stochastic-robust problems (\Cref{mixed_section}). The case study in \Cref{case_study_section} first presents results from tuning the mixed stochastic-robust model (\Cref{weight_section}) and then how it compares with the other described modelling approaches (\Cref{comparison_section}).

\section{Modelling}\label{model_section}
The perspective taken in this paper is that of a system operator aiming at optimally scheduling and balancing a completely renewable system dominated by hydropower. The system is scheduled to be in balance according to the net load forecast in the day-ahead planning stage, and symmetric spinning reserve capacity is procured to ensure the balancing capabilities of the system. The existence of some variable generation components in the system is not modelled explicitly, but manifests as uncertainty in the net load. The forecast errors are seen as the main factors of this uncertainty, and are therefore the drivers behind the need for balancing services. The forecast errors in the net load become known after the scheduling step, and so the operator must use the procured reserve capacity to balance the system in the cheapest way possible. 

\subsection{Deterministic day-ahead scheduling problem}\label{deterministic_section}
The simple deterministic short-term scheduling problem for the system operator,

\begin{equation}\label{da_problem}
\begin{aligned}
&\min_{\vec{x}} Z^{da}(\vec{x})\\
&\vec{x}\in \X, \\
\end{aligned}
\end{equation}

aims to minimize the cost of using water to cover the required net load and spinning reserve requirements while respecting the physical constraints of the system. It is formulated as

\begin{equation}\label{da_obj}
\min \sum_{m\in\M}\left(-WV_{m}v_{m,T+1} +\sum_{t\in\T}F_{t}\left(C^{b}q_{mt}^{b} + C^{o}q_{mt}^{o}\right)\right)
\end{equation}

S.t.
\begin{align}
&q_{mt}^{in} = \sum_{i\in\I{m}{d},n\in\N{i}}q^{d}_{int} + \sum_{j\in\I{m}{b}}q^{b}_{jt} + \sum_{k\in\I{m}{o}}q^{o}_{kt} \label{water_in} \\
&q_{mt}^{out} = \sum_{n\in\N{m}}q^{d}_{mnt} + q^{b}_{mt}+ q^{o}_{mt} \label{water_out}\\
&v_{m0} = V_{m}^{0} \label{init_vol} \\
&v_{m,t+1} - v_{mt} - F_{t}q_{mt}^{in} + F_{t}q_{mt}^{out}  = F_{t}I_{mt}  \label{hydro_bal} \\
&p_{mt} = \sum_{n\in\N{m}}E_{mn}q^{d}_{mnt} \label{hydro_prod}\\
&\sum_{m\in \M}p_{mt} = L_{t}  \label{power_bal}\\
&p_{mt} + r_{mt} \leq P_{m}  \label{prod_upper}\\
&p_{mt} - r_{mt} \geq 0   \label{prod_lower} \\
&\sum_{m\in\M}r_{mt} \geq R_{t} \label{res_req}\\
&v_{mt} \leq V_{m}  \label{max_vol1} \\
&q^{d}_{mnt} \leq Q_{mn}^{d} \label{max_discharge}\\
&q^{b}_{mt} \leq Q_{m}^{b} \label{max_bypass} \\
&q^{o}_{mt} \leq Q_{m}^{o}  \label{max_spill}\\
&p_{mt} \leq P_{m} \label{max_prod} \\
&\text{All variables} \geq 0. \label{da_pos_var}
\end{align}

All symbols in uppercase are input parameters to the model, while lowercase symbols represent the decision variables $\vec{x}$. The model is defined for the hydropower modules $m \in \M$ over the time periods $t\in \T$. The water may be moved between modules through three different waterways: flow through the turbine $q^d$, flow through the bypass gate $q^b$ and spillage $q^o$. The turbine flow is separated into several segments $n\in \N{m}$. The objective (\ref{da_obj}) of the model is to minimize the total cost of using water according to the water value $WV$ and end volume $v_{T+1}$, as well as the penalties for using the bypass and spillage waterways. Constant water values are used in this model, but cutting plane descriptions of the end value of the water may also be used, such as in \cite{Helseth2018a}. \Cref{water_in} sums up the total water entering each module from the set of connected upstream modules $\mathcal{I}$ in each time step, while \cref{water_out} calculates the amount of water released. The mass balance of each reservoir is preserved through \cref{hydro_bal}, where the increase in the volume $v$ over a time step with duration $F$ is equal to the natural inflow $I$ in addition to the net controlled inflow. The starting reservoir content $V^{0}$, set through \cref{init_vol}, serves as the initial condition for the system. The relation between water discharged through the turbine and the power $p$ produced by the generator is modelled as a piece-wise linear constraint in \cref{hydro_prod}, where the efficiency $E$ is decreasing for increasing discharge segment number $n$ to ensure convexity of the problem. The total power balance in \cref{power_bal} ensures that the forecasted net load $L$ is met, and so the system is scheduled to be in balance given the current information. \Cref{prod_upper,prod_lower} bound the available spinning reserve capacity $r$ of the modules based on the maximal production capacity $P$. Enough symmetric spinning reserve capacity must be allocated to satisfy the static reserve requirement $R$ in \cref{res_req}. \Crefrange{max_vol1}{max_prod} are the upper bounds of the variables based on the physical capacities of the hydropower modules, and \cref{da_pos_var} ensures non-negative variables.

\subsection{The balancing problem}\label{bal_section}
Balancing the system in real time after some net load deviation $\Delta$ has occurred is necessary to maintain system stability. The decisions $\vec{x}$ made in the day-ahead planning stage will affect the system's ability to perform the balancing actions, and so the balancing problem

\begin{equation}\label{bal_problem}
\begin{aligned}
&\min_{\vec{y}} Z^{bal}(\vec{y})\\
&\vec{y} \in \Y(\vec{x},\Delta), \\
\end{aligned}
\end{equation}

depends on both $\vec{x}$ and $\Delta$. The formulation is similar to the day-ahead scheduling problem described by \crefrange{da_obj}{da_pos_var} with alterations:

\begin{align}
&\min \sum_{m\in\M}\left(-WV_{m} \overline{v}_{m,T+1} +\sum_{t\in\T}F_{t}\left(C^{b} \overline{q}_{mt}^{b} + C^{o} \overline{q}_{mt}^{o}\right)\right)\nonumber\\
&\hspace{12pt plus 0pt minus 0pt}+ \sum_{t\in\T}\left(C^{+} \overline{s}_{t}^{+}+C^{-} \overline{s}_{t}^{-}\right)\label{bal_obj}
\end{align}

S.t.
\begin{align}
&\overline{q}_{mt}^{in} = \sum_{i\in\I{m}{d},n\in\N{i}} \overline{q}^{d}_{int} + \sum_{j\in\I{m}{b}} \overline{q}^{b}_{jt} + \sum_{k\in\I{m}{o}} \overline{q}^{o}_{kt} \quad(\psi_{mt}^{in})  \label{water_in_bal}\\
& \overline{q}_{mt}^{out} = \sum_{n\in\N{m}} \overline{q}^{d}_{mnt} +  \overline{q}^{b}_{mt}+  \overline{q}^{o}_{mt} \hspace{60pt plus 0pt minus 0pt}(\psi_{mt}^{out}) \label{water_out_bal}\\
& \overline{v}_{m0} = V_{m}^{0} \hspace{153pt plus 0pt minus 0pt} (\nu_{m}) \label{init_vol_bal}\\
&  \overline{v}_{m,t+1} -  \overline{v}_{mt} - F_{t} \overline{q}_{mt}^{in} + F_{t} \overline{q}_{mt}^{out}  = F_{t}I_{mt}  \hspace{19pt plus 0pt minus 0pt} (\mu_{mt}) \label{hydro_bal_bal}\\
&  \overline{p}_{mt} = \sum_{n\in\N{m}}E_{mn} \overline{q}^{d}_{mnt} \hspace{99pt plus 0pt minus 0pt} (\eta_{mt})  \label{hydro_prod_bal}\\
&  \sum_{m\in \M} \overline{p}_{mt} +  \overline{s}^{+}_{t} -  \overline{s}^{-}_{t} = L_{t} + \Delta_{t} \hspace{63pt plus 0pt minus 0pt}(\lambda_{t})\label{power_bal_bal}\\
& \overline{p}_{mt} \leq p_{mt} + r_{mt}   \hspace{121pt plus 0pt minus 0pt}  (\rho_{mt}^{+}) \label{prod_upper_bal}\\
& \overline{p}_{mt} \geq p_{mt} - r_{mt} \hspace{121pt plus 0pt minus 0pt} (\rho_{mt}^{-})    \label{prod_lower_bal}\\
& \overline{v}_{mt} \leq V_{m} \hspace{150pt plus 0pt minus 0pt} (\omega_{mt})  \label{max_vol1_bal}\\
&  \overline{q}^{d}_{mnt} \leq Q_{mn}^{d} \hspace{135pt plus 0pt minus 0pt} (\gamma_{mnt}^{d}) \label{max_discharge_bal}\\
& \overline{q}^{b}_{mt} \leq Q_{m}^{b} \hspace{150pt plus 0pt minus 0pt}(\gamma_{mt}^{b}) \label{max_bypass_bal}\\
&\overline{q}^{o}_{mt} \leq Q_{m}^{o}   \hspace{150pt plus 0pt minus 0pt}(\gamma_{mt}^{o}) \label{max_spill_bal} \\
&\text{All variables} \geq 0. \label{bal_pos_var}
\end{align}

The balancing-stage variables $\vec{y}$ are represented by lower case characters with an overline, and are analogous to the first-stage variables without the overline. The power balance in \cref{power_bal_bal} now includes the net load deviation $\Delta$ and the option to spill power and shed load through the variables $ \overline{s}^{\pm}$ for a penalty cost of $C^{\pm}$. The connection to the first-stage decisions are found in the upper and lower production limits, \cref{prod_upper_bal,prod_lower_bal} respectively, which now constrains the balancing power production to be within the limits set by the scheduled production level $p$ and symmetric spinning reserve capacity $r$. Note that the dual variables for each constraint is included in parentheses.

\subsection{Two-stage problem formulations}\label{two_stage_section}

A two-stage stochastic problem may now be constructed based on the formulations in \Cref{deterministic_section,bal_section}. By providing a set $\SC$ of balancing scenarios with net load deviations $\Delta_{st}$ and probabilities $\pi_s$, a copy of the balancing problem \cref{bal_problem} may be introduced for each scenario. This leads to the classical two-stage stochastic problem formulation \cite{Morales2014}:

\begin{equation}\label{stoch_problem}
\begin{aligned}
&\min_{\vec{x},\vec{y}_s} Z^{da}(\vec{x}) + \sum_{s\in\SC} \pi_s Z^{bal}(\vec{y}_s)\\
&\vec{x}\in \X\\
&\vec{y}_s \in \Y(\vec{x},\Delta_s)\qquad\forall s\in\SC. \\
\end{aligned}
\end{equation}

The robust two-stage counterpart to this stochastic formulation is the tri-level problem

\begin{equation}\label{rob_problem}
\begin{aligned}
&\min_{\vec{x}} Z^{da}(\vec{x}) + \max_{\Delta}\min_{\vec{y}} Z^{bal}(\vec{y})\\
&\vec{x}\in \X \\
&\Delta \in \LO\\
&\vec{y} \in \Y(\vec{x},\Delta),
\end{aligned}
\end{equation}

where the net load deviation $\Delta$ is constrained to be part of the uncertainty set $\LO$. We will use the simple formulation first proposed in \cite{Bertsimas2004} to define:

\begin{equation}\label{uncertainty_set}
\begin{aligned}
\LO  :=  &\Big\{\,\Delta_{t} \mid \Delta_{t} = \Lambda\cdot(u_{t}^{+} - u_{t}^{-})\,; \\
&\, \sum_{t\in\T}(u_{t}^{+} + u_{t}^{-}) \leq \Gamma\,;\,u_{t}^{\pm} \in \{0,1\} \,\Big\},\\
\end{aligned}
\end{equation}

where the parameters $\Lambda$ and $\Gamma$ are the maximum net load deviation and the budget of uncertainty, respectively. The binary variables $u_{t}^{\pm}$ signify if a deviation in positive ($u_{t}^{+}=1$) or negative ($u_{t}^{-}=1$) direction has occurred. The robust problem aims to minimize the cost of the worst case net load deviation contained within the uncertainty set.

The min-max-min formulation of the robust optimization problem in \cref{rob_problem} cannot be solved directly. The column-and-constraint generation (CCG) procedure, first proposed in \cite{Zeng2013,LongZhao2012}, is a popular primal decomposition scheme to remedy this. Other solution techniques such as Benders decomposition (see for instance \cite{Bertsimas2013}) and affine policy approximation \cite{Lorca2016,Lorca2017} are not considered in this work due to lack of efficiency for the problem at hand. CCG first requires the inner minimization problem of \cref{rob_problem} to be transformed to its dual maximization form, so that the maximization steps may be combined:

\begin{equation}\label{second_stage_transformation}
\begin{aligned}
&\max_{\Delta\in \LO}\min_{\vec{y}}Z^{bal}(\vec{y})& \Longleftrightarrow\quad &\max_{\Delta\in \LO,\vec{\phi}} W^{bal}(\vec{x},\Delta,\vec{\phi})\\
&\vec{y} \in \Y(\vec{x},\Delta)&&\vec{\phi} \in \Omega. \\
\end{aligned}
\end{equation}

The vector of dual variables is denoted as $\vec{\phi}$, which are bound by the dual constraints in $\Omega$. The dual variables are listed as Greek lower case letters in parentheses behind \crefrange{bal_obj}{max_spill_bal}. Bi-linear terms $\Delta_t\lambda_t$ appear in the objective function $W^{bal}$ when \cref{power_bal_bal} is dualized. The binary definition of $\Delta_t$ in \cref{uncertainty_set} allows for an exact reformulation of the bi-linear problem to a mixed integer linear program (MILP) by using a "big-M" approach. This is done in for instance \cite{Jiang2012}, though other options are available for solving the problem. An alternating direction method was used in \cite{Lorca2015}, while a cutting plane outer approximation was implemented in \cite{Bertsimas2013}. In this paper the exact MILP reformulation will be used, as it can be directly solved with a standard MILP solver. The CCG technique is based on repeatedly solving the dual form of \cref{second_stage_transformation} for iteratively updated first-stage solutions $x\in\X$. The solution yields the realization of the worst case net load deviation $\Delta_j$, which is iteratively added to the master problem 
\begin{equation}\label{rob_two_stage}
\begin{aligned}
&\min_{\theta,\vec{x},\vec{y}_j} Z^{da}(\vec{x}) + \theta\\
&\vec{x}\in \X\\
&\theta \geq Z^{bal}(\vec{y}_j)\qquad\hspace{6pt plus 0pt minus 0pt}\forall j\in\mathcal{J}\\
&\vec{y}_j \in \Y(\vec{x},\Delta_j)\qquad\forall j\in\mathcal{J}.\\
\end{aligned}
\end{equation}

The set $\mathcal{J}$ represents the number of worst case net load deviation scenarios that have been created, and the auxiliary variable $\theta$ is an outer approximation of \cref{second_stage_transformation}. Solving the master problem results in a new first-stage solution $x$, which is used to solve \cref{second_stage_transformation} in the next iteration. When the current value of $\theta$ and $W^{bal}$ have converged within a specified tolerance, the procedure is complete as the optimal value of \cref{rob_problem} has been found.

\subsection{Mixed stochastic-robust problem}\label{mixed_section}

The unified stochastic-robust optimization first proposed in \cite{Zhao2013} combines the stochastic and robust two-stage problem formulations in \cref{stoch_problem,rob_problem} by introducing a scaling factor  $0\leq \beta \leq 1$:

\begin{align}
&\min_{\vec{x},\vec{y}_s}Z^{da}(\vec{x}) + \beta\sum_{s\in\mathcal{S}}\pi_s Z^{bal}(\vec{y}_s) + (1-\beta)\max_{\Delta}\min_{\vec{y}}Z^{bal}(\vec{y})\nonumber\\
&\vec{x}\in \X \nonumber\\
&\vec{y}_s \in \Y(\vec{x},\Delta_s)\qquad\forall s\in\SC\label{unified_problem}\\
&\Delta\in\LO\nonumber\\
&\vec{y}\in\Y(\vec{x},\Delta).\nonumber
\end{align}

The scaling factor is an importance weighting of the worst-case and expected costs, where choosing the edge cases $\beta = 1$ or $\beta = 0$ results in pure stochastic or pure robust objective functions, respectively. The formulation given in \cref{unified_problem} can be regarded as an augmented version of the original robust tri-level problem formulation in \cref{rob_problem}, where the first-stage problem is extended to include the balancing scenarios in $\SC$. The CCG solution techniques can therefore still be used to solve the unified problem, as it is still the min-max-min form. A drawback of the unified problem formulation is the increased size of the master problem in the CCG algorithm, as all constraints related to the exogenously created balancing scenarios are added. This could result in poor computational performance and may call for a second decomposition scheme to solve the master problem itself. 

We propose a similar but novel mixed stochastic-robust formulation that compartmentalizes the model complexity, facilitates reusability and strengthens the robustness of the solution. The complexity of the model is reduced compared to \cref{unified_problem} by separating out the calculation of the robust elements. This is achieved by viewing the CCG algorithm as a robust scenario generator. Solving the pure robust model described in \cref{rob_problem} with the CCG algorithm creates a set $\mathcal{J}$ of robust scenarios. These scenarios are then added to the stochastic model formulation in \cref{stoch_problem}, again scaled with $\beta$:

\begin{align}
&\min_{\vec{x},\vec{y}_s,\vec{y}_j}Z^{da}(\vec{x}) + \beta\sum_{s\in\mathcal{S}}\pi_s Z^{bal}(\vec{y}_s) + (1-\beta)\sum_{j\in\mathcal{J}}\pi_{j}Z^{bal}(\vec{y}_{j})\nonumber\\
&\vec{x}\in \X\nonumber\\
&\vec{y}_s \in \Y(\vec{x},\Delta_s) \qquad\forall s\in\SC\label{mixed_problem}\\
&\vec{y}_{j} \in \Y(\vec{x},\Delta_{j})\qquad\forall j\in\mathcal{J}.\nonumber
\end{align}

The robust scenarios are assumed to be equiprobable, $\pi_j = 1/|\mathcal{J}|$. The separate generation of the robust scenarios simplifies the solution strategy of the new mixed stochastic-robust model in \cref{mixed_problem}, as this extensive form representation can be solved directly by a linear programming solver. The separation may lead to lower calculation times, though it will be highly dependent on the individual case. The construction of the robust scenarios no longer depends on the exogenously generated scenarios, making them potentially valid for other days with similar net load profiles. The set of scenarios $\mathcal{J}$ can then be used as an initial set of constraints in the solution of the pure robust model. This sharing of contingency events between time periods has also been proposed for the long-term hydrothermal planning model in \cite{Street2017}, which incorporates CCG in a SDDP framework. 

The mixed model in \cref{mixed_problem} minimizes the expected cost of all robust scenarios instead of minimizing the cost of the worst case scenario. This will increase the cost of the very worst robust scenario, but also force all robust scenarios to make rational use of their water. The balancing problem constraints \crefrange{water_in_bal}{bal_pos_var}, constitute a weak coupling to the first-stage, in the sense that the problem has complete recourse. It is possible to shed load and pass water through unfavourable spill and bypass gates to preserve feasibility. When only the maximal robust scenario cost is minimized, the weak coupling causes most of the robust scenarios to utilize these shortcuts to some extent. Enforcing the minimization of the expected cost increases the robustness of the model, while the overall conservativeness can still be regulated through the weight $\beta$. The coupled model proposed in \cite{Liu2017} also minimizes the expected value of the robust scenarios, but in their case the robust scenarios are iteratively added to the set of exogenously generated scenarios through the CCG procedure. The approach still has the potential tractability issues of the original unified problem formulation, especially since the convergence of the CCG algorithm is unproven in their proposed framework.


\section{Case study}\label{case_study_section}

The focus of this case study is on the quality of the solution of the mixed stochastic-robust model proposed in \Cref{mixed_section} in terms of costs and geographical coordination of the procured reserves, and how this compares to the other proposed models in \Cref{model_section}. All optimization models have been implemented in the Pyomo modelling package for Python \cite{Hart2011,Hart2017} using the MILP solver CPLEX 12.8 \cite{cplex}.

\begin{figure}[!t]
    \centering 
    \includegraphics[width=0.8\columnwidth]{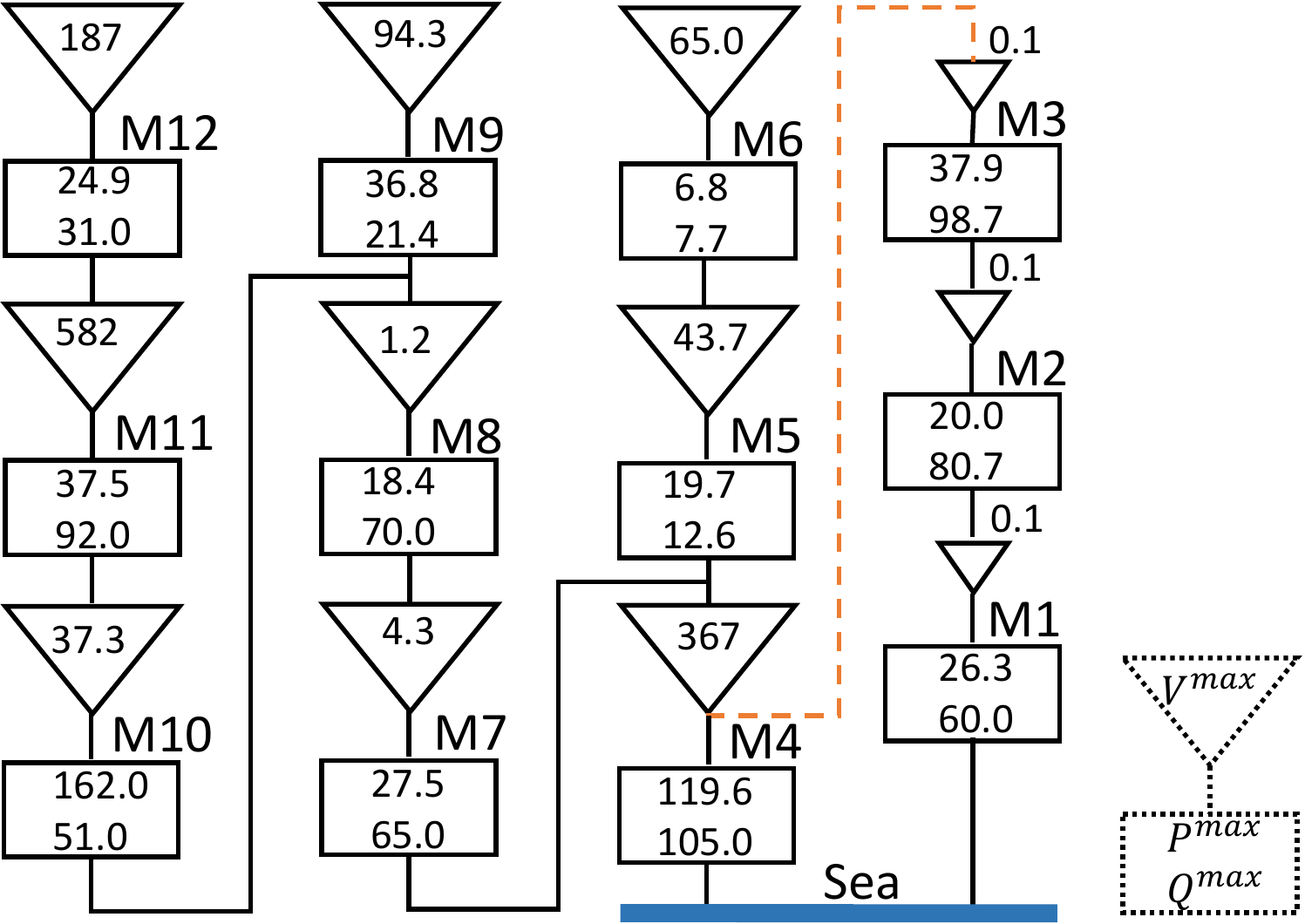}
\caption{Sketch of the hydropower topology. Reservoirs (triangles), power plants (rectangles), and water routes for discharge (solid lines) and bypass (dashed line) are shown together with maximal values for discharge ($\text{m}^3$/s), production (MW) and reservoir volumes (M$\text{m}^3$).}
\label{fig:topology}
\end{figure}

\begin{figure}[!t]
    \centering
    \includegraphics[width=\columnwidth]{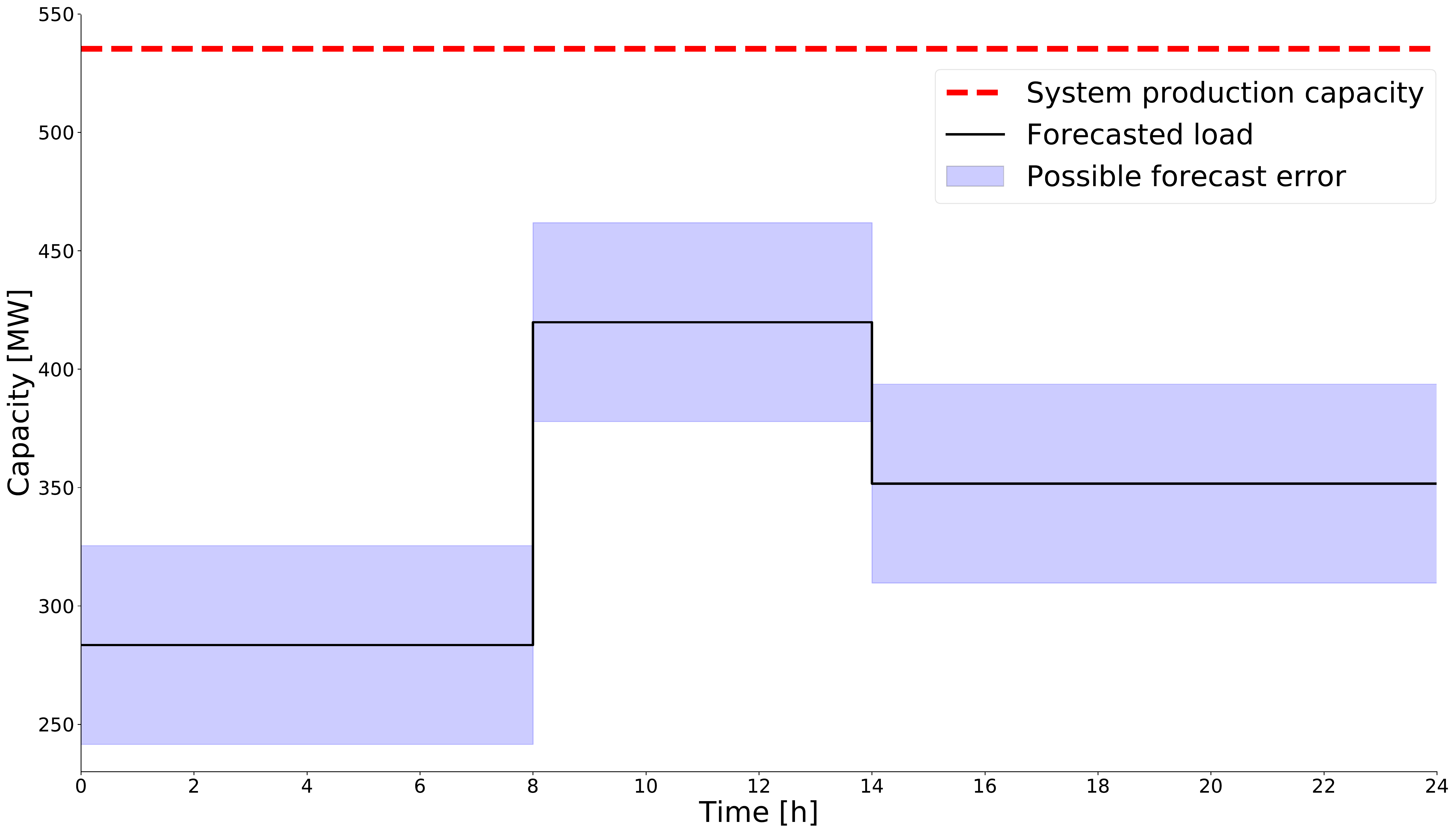}
\caption{Forecasted system net load with the region of possible deviations and the total installed capacity of the system.}
\label{fig:load}
\end{figure}

The hydropower system used in the study is shown in \Cref{fig:topology}. It is based on a real watercourse in Norway, and consists of 12 modules with a total installed production capacity of $537.4$ MW. The initial reservoir volume of every module is set to $65\%$ of its maximal storage capacity, which represents a normal hydrological situation during the winter in Norway. Water values are calculated by a long-term hydropower scheduling model described in \cite{Helseth2018a}, and are measured in arbitrary monetary units (mu) per M$\text{m}^{3}$ in the range of 1200 - 9000 mu/M$\text{m}^{3}$. The penalty for shedding load and spilling power is chosen to be 3000 mu/MW and 1000 mu/MW, respectively. The time horizon is set to 24h with hourly resolution, and the forecasted net load profile is shown in \Cref{fig:load}. The parameters of the robust uncertainty set $\LO$ are chosen to be $\Lambda = 42$ MW and $\Gamma = 6$, resulting in a worst case balancing situation with up to 6 hours with net load deviations of $42$ MW. This maximal net load deviation is $10\%$ of the 420 MW peak forecasted net load, and the duration of the forecasted net load peak corresponds to the budget of uncertainty $\Gamma$.

There are two measures that are used to quantify the quality of the proposed coordination of procured reserves between the hydropower modules; The cost of procuring the reserves $K$ and the following cost of balancing the system $B$. The cost of procuring reserves is defined as the increase in the first-stage objective function $Z^{da}$ relative to the cost of the deterministic model in \cref{da_problem} solved without any reserve requirement, $Z^{da}_{0}$:

$$K = Z^{da} - Z^{da}_{0}.$$

This cost represents the opportunity cost of procuring the reserves, and should be paid as compensation to the owners of the hydropower plants for following the shifted schedule. For a given net load deviation realization $\Delta_{it}$, the cost $B_i$ of balancing this deviation is calculated by solving the primal balancing problem in \cref{bal_problem} for the given production set points $p$ and allocated reserves $r$. This yields the objective function $Z^{bal}_i$, which is normalized by the objective function given perfect foresight, $Z^{PF}_{i}$, to produce the balancing cost:

$$B_{i} = Z^{bal}_{i} - Z^{PF}_{i}.$$

The perfect foresight cost is found by relaxing the production limit constraints of the balancing problem, \cref{prod_upper_bal,prod_lower_bal}, to let every module produce between zero and maximum capacity. The sum of the procurement cost and the subsequent balancing cost is the normalized total system cost,

$$U_i = K + B_i.$$

$K$ is easily found by direct calculation, whereas the balancing costs $B_i$ must be estimated by simulation. In this case study, 1000 different sampled net load deviations were used to measure the balancing cost of the given schedule. The sampled net load deviations were drawn from a normal distribution with mean zero and standard deviation $\Lambda/2.5$, where any values $|\Delta_{it}| \geq \Lambda$ were truncated. This was done to ensure that only the quality of the coordination of the reserve procurement was measured, and not simply the amount of reserves procured. The exogenously generated scenarios used in the stochastic, unified and mixed stochastic-robust models were constructed from the same truncated normal distribution, but redrawn based on a different seed for the random number generator. The same exogenous scenarios were used in all of these models.

\subsection{Choosing \texorpdfstring{$\beta$}{b} in the mixed stochastic-robust model}\label{weight_section}
The pure robust model in \cref{rob_problem} was solved to an absolute convergence tolerance of 1 mu, with an absolute MIP gap of 0 and integer tolerance of $10^{-9}$ for the second-stage problem. A total of 23 iterations of the CCG algorithm were performed, resulting in 23 robust scenarios. An additional 50 equiprobable scenarios were generated to form the set $\SC$ of balancing scenarios. To determine the optimal level of $\beta$ in the mixed model, 11 simulation runs with values of $\beta$ ranging from 0 to 1 in steps of 0.1 were performed. As described previously, each simulation run consisted of calculating $B_i$ for 1000 sampled balancing scenarios. To test the flexibility and robustness of the model to uncertainty in the underlying net load deviation distribution, another 11 simulation runs were performed with scenarios drawn from a uniform distribution, $-\Lambda \leq \Delta_{it} \leq \Lambda$. The numerical data is presented in \Cref{table:beta_sim}. 


\begin{table}[!t]
\centering
\caption{Costs, measured in monetary units, of the simulation runs with different weights $\beta$. The 1000 simulation scenarios were drawn from a normal distribution in (a) and a uniform distribution in (b).}\label{table:beta_sim}
\subfloat[]{
\begin{tabular}{SSSSSS} \toprule
    {$\beta$} & {$K$} & {$U^{max}$} & {$U^{mean}$} & {$U^{med}$} & {$\sigma(U_i)$}\\ \midrule
    0.0& 60.7&    84.8&  64.6&  63.4 & 3.8 \\
    0.1& 54.1&    75.5&  57.6&  56.6 & 3.6 \\ 
    0.2& 39.2&    74.0&  43.2&  41.7 &  4.5\\
    0.3& 23.8&    67.5&  28.9&  27.1 &  5.5\\  
    0.4& 18.4&    62.1&  23.7&  21.9 &  5.8\\
    0.5& 16.1&    59.8&  21.9&  20.0 &  6.0\\
    0.6& 12.9&    60.5&  19.6&  17.5 &  6.7\\
    0.7& 11.0&    60.2&  18.2&  16.0 &  7.1\\
    0.8& 9.6&    59.6&  17.7&  15.3  & 7.6\\
    0.9& 8.9&    59.4&  17.6&  15.2  & 7.9\\
    1.0& 7.4&    78.8&  18.8&  15.3  & 10.8\\ \bottomrule
\end{tabular}}
\hfill
\subfloat[]{
\begin{tabular}{SSSSSS} \toprule
    {$\beta$}& {$K$} & {$U^{max}$} & {$U^{mean}$} & {$U^{med}$} & {$\sigma(U_i)$} \\ \midrule
    0.0& 60.7&   96.9&  73.2&  72.2 & 6.2  \\
    0.1&  54.1&  96.5&  66.0&  64.7 & 6.7  \\ 
    0.2&  39.2&  85.3&  53.5&  51.8 & 8.2  \\
    0.3&  23.8&  91.2&  41.6&  39.6 & 10.3  \\  
    0.4&  18.4&  85.1&  37.7&  36.1 & 10.7  \\
    0.5&  16.1&  82.9&  37.2&  35.6 & 11.1  \\
    0.6&  12.9&  83.1&  36.7&  35.1 & 12.1  \\
    0.7&  11.0&  85.1&  36.8&  35.1 & 12.7  \\
    0.8&  9.6&   86.5&  37.8&  35.9 & 13.4  \\
    0.9&  8.9&  87.7&  38.9&  37.0  & 13.9 \\
    1.0&   7.4& 105.4& 45.1&  42.8  & 17.3 \\ \bottomrule
\end{tabular}}
\end{table}


It is clear that a lower value of $\beta$ gives a lower sample standard deviation $\sigma$ at the expense of a higher base procurement cost $K$. The model has the lowest procurement cost $K$ when $\beta=1$, however, there is a noticeable increase in the maximal cost and standard deviation for this value. This is even more apparent in the case with uniformly distributed balancing scenarios. Based on the median and average values of the data, the optimal choice is $\beta = 0.9$ in the case of normally distributed net load deviations, and $\beta=0.6$ when the simulated scenarios are drawn from a uniform distribution. The fact tat a non-zero robust weight minimizes the total costs in both cases signifies a benefit of strengthening the stochastic model with robust scenarios. This is especially true when the probability distribution used to generate the scenarios for the stochastic model is inaccurate, which can be the case when the underlying data is of poor quality. The cost of increasing the robustness of the model is also moderate in the region $0.7\leq\beta\leq1$.

\subsection{Model comparison}\label{comparison_section}
\begin{figure}[!t]
    \includegraphics[width=\columnwidth]{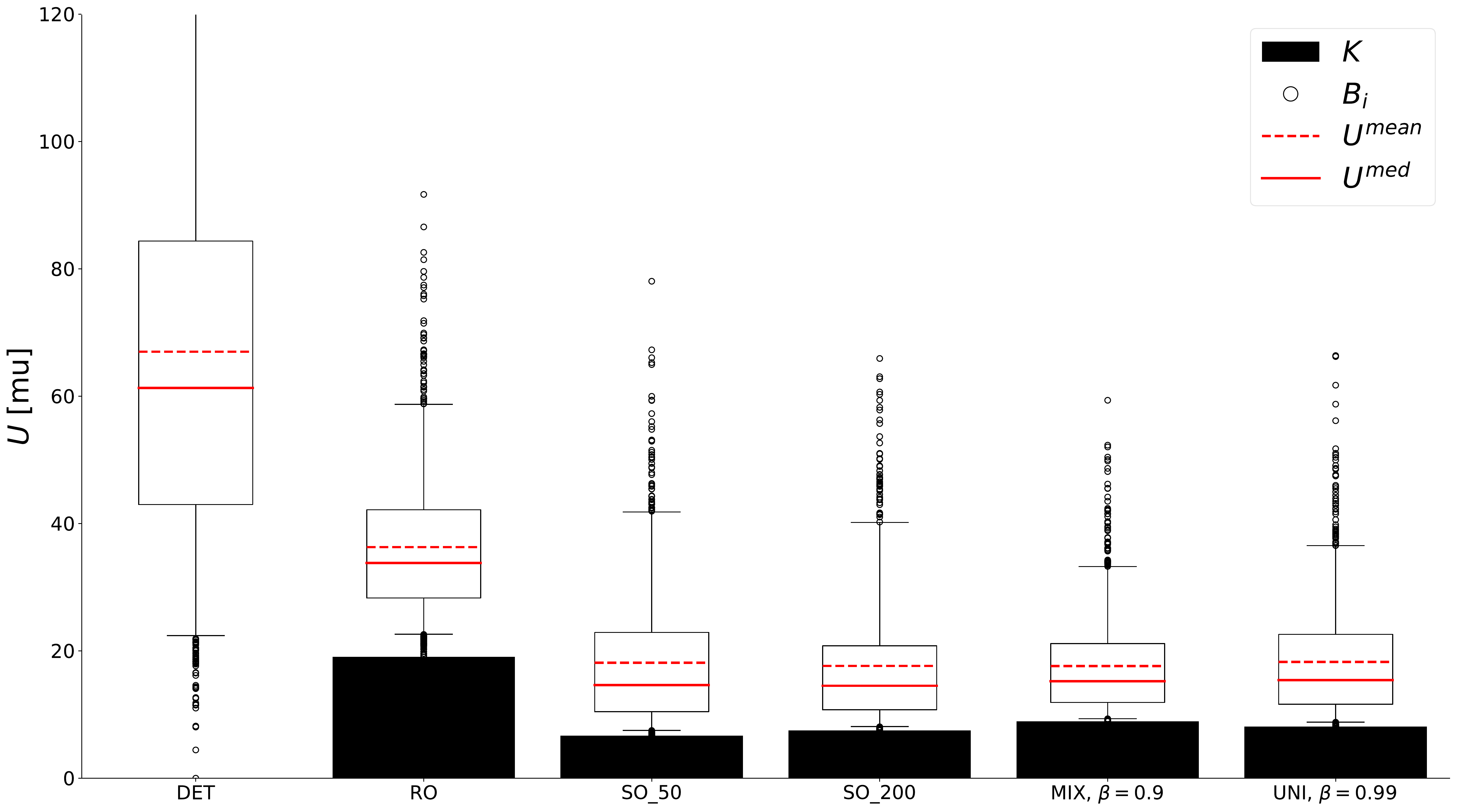}
\caption{A box plot of the simulated balancing costs $B_i$ on top of the static procurement costs $K$ for 1000 simulated scenarios drawn from a normal distribution. The whiskers of the boxes are at the $5\%$ and $95\%$ percentiles.}
\label{fig:box_compare}
\end{figure}

To compare the mixed stochastic-robust model with the deterministic (\cref{da_problem}), pure robust (\cref{rob_problem}), pure stochastic (\cref{stoch_problem}) and unified stochastic-robust (\cref{unified_problem}) models, the same simulation run of 1000 balancing scenarios drawn from the normal distribution was used to simulate the costs. The deterministic model was solved with a static reserve requirement of $R_t = \Lambda = 42$MW and the pure robust model was solved with the same tolerances described in \Cref{weight_section}. The pure stochastic model was solved first with 50 and then with 200 scenarios, and the unified model used 50 scenarios. Note that the weight $\beta=0.99$ was chosen for the unified model after performing an analysis similar to the description in \Cref{weight_section} with a resolution of 0.01 in the region $0.9\leq\beta\leq1.0$. The results of the simulation are shown in \Cref{fig:box_compare} as a box plot, while numerical details are given in \Cref{table:model_sim}. 

\begin{table}[!t]
\centering
\caption{Costs, measured in monetary units, of the simulation runs based on 1000 sampled scenarios from a normal distribution.}\label{table:model_sim}
\begin{tabular}{SSSSSS} \toprule
    {\text{Model}} & {$K$} & {$U^{max}$} & {$U^{mean}$} & {$U^{med}$} & {$\sigma(U_i)$}\\ \midrule
    \text{Deterministic}& 0 &    210.3&  67.0&  61.3 & 33.9 \\
    \text{Robust}&  19.1 &  91.7&  36.3&  33.8 & 11.4\\ 
    \text{Stochastic, 50 scen.}& 6.7 &   78.1&  18.1&  14.6  & 10.8\\
    \text{Stochastic, 200 scen.}& 7.5 &   65.9&  17.6&  14.5 &  10.1\\
    \text{Mixed, $\beta=0.90$}& 8.9 &  59.4&  17.6&  15.2  & 7.9\\  
    \text{Unified, $\beta=0.99$}& 8.1 &  66.4&  18.3&  15.4 &  9.1\\  \bottomrule
\end{tabular}

\end{table}

\begin{figure}[!t]
\subfloat[Deterministic model]{\includegraphics[width=\columnwidth]{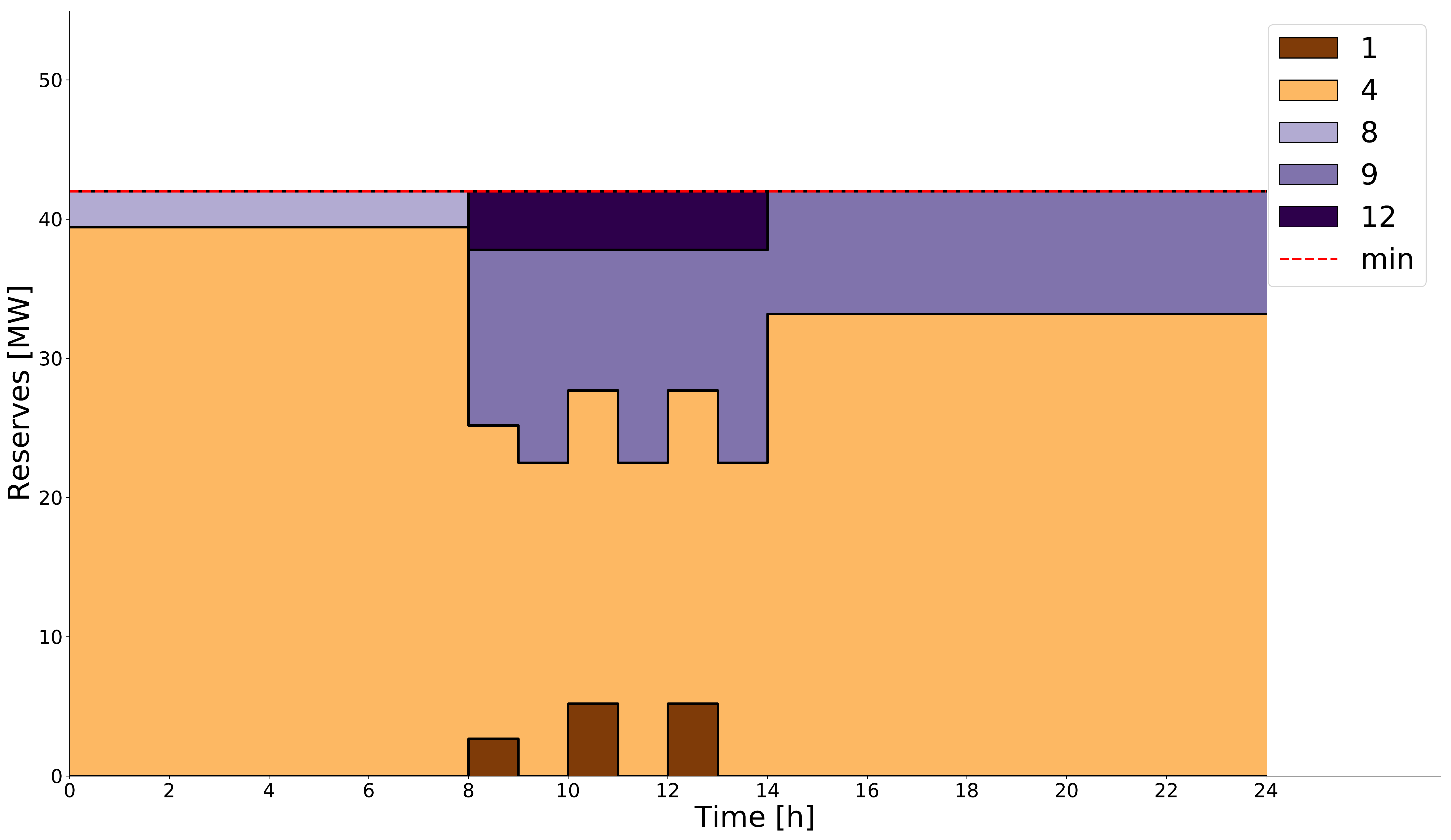}\label{fig:reserves_det}}\\
\subfloat[Mixed stochastic-robust model, $\beta=0.9$]{\includegraphics[width=\columnwidth]{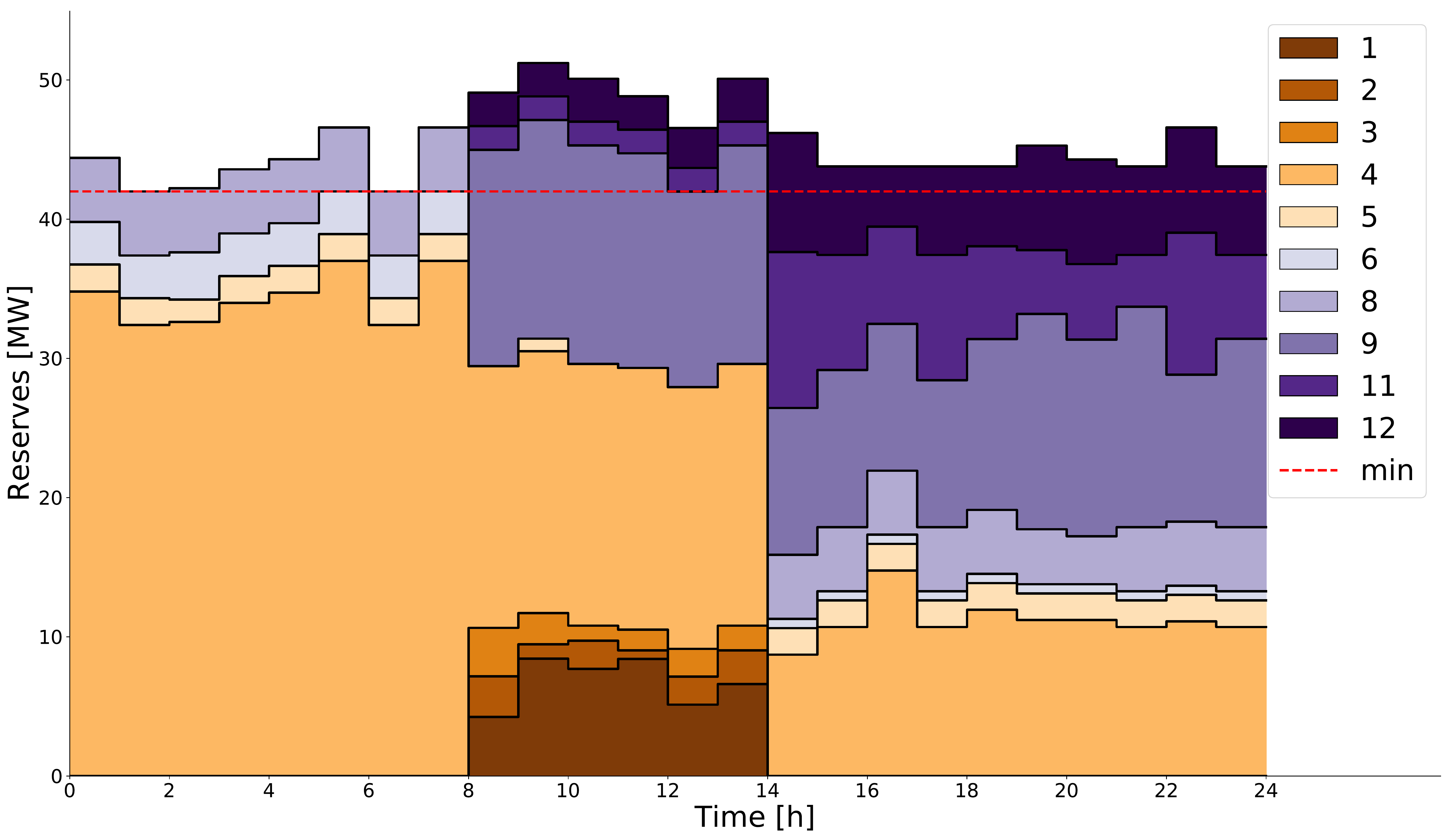}\label{fig:reserves_sro}}
\caption{Reserve allocation on the 12 hydropower modules for two different models. The minimal reserve requirement is shown as a red dashed line.}\label{fig:reserves}
\end{figure}

The pure stochastic models have slightly lower median costs compared to the mixed model, but both the standard deviation and average cost is higher. The spread of the costs in the two stochastic models are very similar, which shows how hard it can be to increase the robustness of a stochastic model by only adding more scenarios. The decrease in the spread of the simulated costs found in the mixed model is therefore valuable. The pure robust model has a spread similar to those of the pure stochastic models, but has significantly higher procurement costs. This shows its overly conservative nature. The unified stochastic-robust model has an optimal value of $\beta$ which is heavily skewed towards the stochastic objective, which was only revealed by simulating the region $0.9\leq\beta\leq1.0$ with a step size of 0.01. This still yields a poorer solution than the mixed model, and arguably the stochastic models. All of the two-stage models handily outperform the deterministic model, which is 3.8 times more expensive on average compared to the mixed model. The procurement cost is zero in the deterministic case, meaning that there are enough spinning reserves in the system to cover the static reserve requirement of $42$ MW without changing the set points of the generators. However, the lack of coordination between modules when reserves are allocated ensures that the deterministic model often encounters a relatively costly balancing step. 

To better understand the importance of coordinating reserve procurement, \Cref{fig:reserves} shows the distribution of the allocated reserves among the 12 hydropower modules in the deterministic and mixed models. The deterministic model relies heavily on module 4 throughout the day, and only procures reserves for four other modules. The mixed model also depends on module 4 for the first eight hours, but other modules further up in the system alleviate the pressure on this module for the remainder of the day. Modules 1, 2 and 3 are connected behind a bypass gate in a string with limited storage capacity in between. These modules contribute to the reserve pool during the peak hours, and are coordinated to ensure that most of the water can be passed through the entire string if activated. Modules 5 and 6 also contribute in the mixed model solution, and allows for some refilling of water in module 4. The mixed model procures extra reserves in the peak hours due to its robust influence. Net load deviations in these hours are clearly important when protecting against a worst case solution, and having 10 MW extra reserves at hand gives more flexibility when deciding which modules should ramp up or down in the different balancing scenarios.

\section{Conclusion}\label{conclusion_section}
In this paper we have shown that the proposed mixed stochastic-robust optimization model can outperform the deterministic, stochastic, robust and unified stochastic-robust models for the hydropower reserve procurement and power scheduling problem. The simplicity of the mixed model is appealing for several reasons, such as ease of implementation and reusability of robust scenarios. Potential for decreased solution time is present, though this will depend heavily on the specific case and convergence criterion used for the pure robust model. Furthermore, the effect of coordinating reserve capacity among the modules in a cascaded hydropower system has been demonstrated. The loss of the connection to the physical system that occurs when activation is not considered leads to poor reserve allocation in terms of balancing costs. 

The case study presented in this paper is based on a moderate case regarding the initial state and energy available
in the system. Solving the daily scheduling problem over a longer period of time, such as with a rolling horizon simulator, will give a more complete picture of the effects of spatially coordinating the reserve capacity allocation. Investigating the consequence and response of a grid in a multi-area model is also an interesting future avenue of research.


%




\ifCLASSOPTIONcaptionsoff
  \newpage
\fi



%
\bibliographystyle{IEEEtran}
\bibliography{refs}

%







\end{document}